\documentclass[a4paper, 12pt]{amsart}

\usepackage[T1]{fontenc}
\usepackage[utf8]{inputenc}
\usepackage[english]{babel}

\usepackage{amssymb,amsmath,amsfonts,amsthm}

\usepackage{xcolor}
\usepackage[normalem]{ulem} 

\usepackage{tikz}
\usetikzlibrary{matrix,arrows,decorations.pathmorphing}
\usepackage{tikz-cd}

\usepackage[mathscr]{euscript} 

\usepackage{enumitem} 

\theoremstyle{plain} 
\newtheorem{thm}{Theorem}[section]
\newtheorem{prop}[thm]{Proposition}
\newtheorem{lmm}[thm]{Lemma}
\newtheorem{cor}[thm]{Corollary}
\newtheorem*{mainthm}{Main Theorem}

\theoremstyle{definition} 
\newtheorem{dfn}[thm]{Definition}

\theoremstyle{remark} 
\newtheorem{rmk}[thm]{Remark}
\newtheorem*{notat}{Notation}

\DeclareMathOperator{\sdiv}{div}
\DeclareMathOperator{\car}{char}
\DeclareMathOperator{\Ext}{Ext}
\DeclareMathOperator{\h}{H}
\DeclareMathOperator{\Hilb}{Hilb}

\DeclareMathOperator{\id}{id}

\DeclareMathOperator{\Quot}{Quot}
\DeclareMathOperator{\rk}{rk}
\DeclareMathOperator{\Spec}{Spec}
\DeclareMathOperator{\Supp}{Supp}
\DeclareMathOperator{\Sym}{Sym}

\newcommand{\mcD}{\mathscr{D}} 
\newcommand{\mcE}{\mathscr{E}}                
\newcommand{\mcF}{\mathscr{F}} 
\newcommand{\mcG}{\mathscr{G}}
\newcommand{\mcH}{\mathscr{H}}

\newcommand{\mcK}{\mathscr{K}}
\newcommand{\Kvar}{\mathrm{K}_0(\mathsf{Var})}
\newcommand{\mcL}{\mathscr{L}}
\newcommand{\mcO}{\mathscr{O}}
\newcommand{\mcQ}{\mathscr{Q}}

\newcommand{\N}{\mathbb{N}}
\newcommand{\Z}{\mathbb{Z}}
\newcommand{\Q}{\mathbb{Q}}

\newcommand{\pr}{\mathbb{P}}

\renewcommand{\epsilon}{\varepsilon}

\renewcommand{\phi}{\varphi}

\newcommand{\mb}[1]{\mathbf{#1}}
\newcommand{\mr}[1]{\mathrm{#1}}

\newcommand{\iso}{\xrightarrow{{}_\thicksim}}
\newcommand{\into}{\hookrightarrow}

\newcommand{\xto}{\xrightarrow}

\newcommand{\msout}{\bgroup\markoverwith{\textcolor{blue}{\rule[0.55ex]{2pt}{0.4pt}}}\ULon}
\newcommand{\mcross}{\bgroup\markoverwith{\textcolor{blue}{$\times$}}\ULon}

\newcommand{\bsout}{\bgroup\markoverwith{\textcolor{green}{\rule[0.55ex]{2pt}{0.4pt}}}\ULon} 

\newcommand{\bcross}{\bgroup\markoverwith{\textcolor{green}{$\times$}}\ULon} 

\newcommand{\fsout}{\bgroup\markoverwith{\textcolor{red}{\rule[0.55ex]{2pt}{0.4pt}}}\ULon} 

\newcommand{\fcross}{\bgroup\markoverwith{\textcolor{red}{$\times$}}\ULon} 

\begin{document}

\title[Quot schemes of zero-dimensional quotients on a curve]{On the motive of Quot schemes of zero-dimensional quotients on a curve}

\author{Massimo Bagnarol}
\address[M.\ Bagnarol]{SISSA -- International School for Advanced Studies, Via Bonomea 265, 34136 Trieste, Italy}
\email{mbagnarol@sissa.it}

\author{Barbara Fantechi}
\address[B.\ Fantechi]{SISSA -- International School for Advanced Studies, via Bonomea 265, 34136 Trieste, Italy}
\email{fantechi@sissa.it}

\author{Fabio Perroni}
\address[F.\ Perroni]{Dipartimento di Matematica e Geoscienze, Università degli Studi di Trieste, via Valerio 12/1, 34127 Trieste, Italy}
\email{fperroni@units.it}


\keywords{Quot schemes; Grothendieck ring of varieties.}

\begin{abstract}
For any locally free coherent sheaf on a fixed smooth projective curve, 
we study the class, in the Grothendieck ring of varieties, of the Quot scheme that parametrizes
zero-dimensional quotients of the sheaf. We  prove that this class
depends only on the rank of the  sheaf and on the length of the quotients.
As an application, we obtain an explicit formula that expresses it in terms of the symmetric products of the curve.
\end{abstract}

\maketitle

\tableofcontents


\section{Introduction}
Let $C$ be a smooth projective curve over an algebraically closed ground field $k$, and let $\mcE$ be a locally free sheaf of rank $r$ on $C$. 
For any $n\in \Z_{>0}$, let  $\Quot^n_{C/k}(\mcE)$ be the Quot scheme which parametrizes coherent quotients of $\mcE$ with finite support and $n$-dimensional space of global sections. 
This Quot scheme is a smooth projective $k$-variety, and it follows easily from its definition that there is a natural isomorphism $\Quot^n_{C/k}(\mcE) \cong \Quot^n_{C/k}(\mcE \otimes \mcL)$ for any invertible sheaf $\mcL$ on $C$; in particular, if $r=1$, $\Quot^n_{C/k}(\mcE)\cong \Quot^n_{C/k}(\mcO_C)$. 
In the case where $r> 1$, the isomorphism class of $\Quot^n_{C/k}(\mcE)$ depends on $\mcE$, as one already sees when  $n=1$, in which case  $\Quot^1_{C/k}(\mcE) \cong \pr(\mcE)$, the projective space bundle associated to $\mcE$ (see Section \ref{n=1}).
Our aims are to study the class of $\Quot^n_{C/k}(\mcE)$ (where $r\geq 2$) in the Grothendieck ring of $k$-varieties and to compute it in terms of the classes $[\Sym^{m}(C)]$, $m\geq 0$.
In particular, we prove the following result.
\begin{mainthm}
Under the previous hypotheses, the equality 
\[
[\Quot^n_{C/k}(\mcE)] = [\Quot^n_{C/k}(\mcO_C^{\oplus r})] 
\] 
holds true in the Grothendieck ring $\Kvar$ of $k$-varieties. 
\end{mainthm}

\subsubsection*{Organization of the paper} 
In Section \ref{notations} we present some well-known facts about Quot schemes on smooth projective curves. In particular, in Section \ref{sigma} we recall the existence of a natural morphism $\sigma: \Quot^n_{C/k}(\mcE) \to \Sym^n(C)$, which we describe explicitly in Section \ref{explsigma}. Finally, in Section \ref{n=1} we consider $\Quot^1_{C/k}(\mcE)$, and we show that, in general, it is not isomorphic to $\Quot^1_{C/k}(\mcO_C^{\oplus r})$.

A more detailed study of the morphism $\sigma$ is the subject of Section \ref{sigmafiber}, where we show that the fibers of $\sigma$ only depend (up to isomorphism) on the rank of the locally free sheaf $\mcE$.

Section \ref{coresect} contains the proof of our main theorem. As an application, we explicitly compute $[\Quot^n_{C/k}(\mcE)] \in \Kvar$ and we prove that $\Quot^n_{C/k}(\mcE)$ is irreducible.

\subsubsection*{Acknowledgements} 
The research was partially supported  by the national project
PRIN 2015EYPTSB-PE1 "Geometria delle variet\`a algebriche", by
FRA 2018 of the University of Trieste, and by the group GNSAGA of INDAM. 
Part of the paper was written while the last author was visiting  Fudan University (Shanghai) and he would like to thank Prof.\ Meng Chen for the invitation and hospitality.


\section{Notations and basic results} \label{notations}
In this section, we recall some basic results that are relevant for us and we fix the notation.
For the proofs and for more details we refer to \cite{Gro61}.
Throughout the paper, we work over an algebraically closed ground field $k$. By a variety we mean a reduced separated scheme of finite type over $k$, not necessarily irreducible.

Let $X$ be a projective scheme. Let $\mcO_X(1)$ be a very ample line bundle on $X$,
let $P\in \Q[t]$ be a polynomial with rational coefficients, and let $\mcF$ be a coherent sheaf on $X$.
We denote by $\Quot^P_{X/k}(\mcF)$ the Quot scheme that parametrizes coherent quotients of $\mcF$ with Hilbert polynomial $P$. Let us recall that $\Quot^P_{X/k}(\mcF)$ is a projective scheme, which represents the contravariant functor that associates to any locally noetherian scheme $S$ the set of isomorphism classes of $S$-flat coherent quotients $q: \mcF_S \to \mcH$, such that the Hilbert polynomial of $\mcH_s$ is equal to $P$, for all $s\in S$.
Here, $\mcF_S$ (respectively $\mcH_s$) is the pullback of $\mcF$ to $S\times_k X$ under the projection onto the second factor (respectively the pullback of $\mcH$ to $X_s$). 
In particular, the identity morphism of $\Quot^P_{X/k}(\mcF)$ corresponds to the universal quotient $u: \mcF_{\Quot^P_{X/k}(\mcF)} \to \mcQ$.

A similar result holds true if $X$ is replaced by a quasi-projective scheme $U$. 
In this case, one defines a functor as before, with the additional requirement that $\mcH$ has proper support over $S$. 
Then, this functor is representable by a quasi-projective scheme $\Quot^P_{U/k}(\mcF)$. The relation between the two constructions is given by the following result (see also \cite{Nit05}).

\begin{thm}\label{GN}
Let $X$ and $\mcF$ be as before, and let $U\subseteq X$ be an open subscheme. Then $\Quot^P_{U/k}(\mcF|_U)$ is naturally an open subscheme of $\Quot^P_{X/k}(\mcF)$.
\end{thm}

Now, let us consider the case where $P$ is a constant polynomial equal to $n\in \Z_{>0}$. 
Then $\Quot^n_{X/k}(\mcF)$ parametrizes coherent quotients of $\mcF$ with finite support and such that the dimension of the space of sections is equal to $n$. Therefore, $\Quot^n_{X/k}(\mcF)$ is independent of $\mcO_X(1)$.

In this paper, the main object of study is $\Quot^n_{C/k}(\mcE)$, where $C$ is a smooth projective curve and $\mcE$ is locally free. 

\begin{notat}
Throughout the article, $\Quot^n_{C/k}(\mcE)$ will be denoted by $\mr{Q}^n_C(\mcE)$. Accordingly, the corresponding universal quotient will be denoted by $u: \mcE_{\mr{Q}^n_C(\mcE)} \to \mcQ$. Whenever $U \subseteq X$ is an open subscheme, we will write $\mr{Q}^n_U(\mcE)$ for the Quot scheme $\Quot^n_{U/k}(\mcE|_U)$, and $u_0: \mcE_{\mr{Q}^n_U(\mcE)} \to \mcQ_0$ for its universal quotient.
\end{notat}

Let us first recall the following fact.

\begin{lmm} \label{smooth}
Let $C$ be a smooth projective curve, and let $\mcE$ be a locally free coherent sheaf of rank $r$ on $C$.
Then $\mr{Q}^n_{C}(\mcE)$ is a smooth variety of dimension $nr$. 
\end{lmm}

\proof
Let  $[q: \mcE \to \mcH]$ be a $k$-rational point of $\mr{Q}^n_C(\mcE)$.
Since the support of $\mcH$ is $0$-dimensional, we have that
\[
\Ext^1(\ker(q),\mcH) \cong \h^1(C,\ker(q)^\vee \otimes \mcH) = 0 \,.
\]
The smoothness now follows from \cite[Prop.\ 2.2.8]{HL10}. Moreover, the dimension coincides with that
of the Zariski tangent space at the point $[q: \mcE \to \mcH]$, which is equal to  $\dim \h^0(C,\ker(q)^\vee \otimes \mcH) =nr$.
\endproof



As a consequence of Theorem \ref{GN} and Lemma \ref{smooth}, $\mr{Q}^n_U(\mcE)$ is a smooth quasi-projective variety, for any $U \subseteq C$ open.

\begin{rmk}
We will prove in Corollary \ref{irreducibility} that, under the above hypotheses, $\mr{Q}^n_C(\mcE)$ is irreducible.
\end{rmk}


\subsection{The morphism $\sigma$}\label{sigma}
In the proof of the Main Theorem we will use the morphism $\frak{N}_{X/k}$ defined in \cite[\S 6]{Gro61},
which will be denoted $\sigma$ in this article. The following result is a special case of Grothendieck's construction.

\begin{prop} 
Let $\mcF$ be a coherent sheaf on $C$. 
Then there exists a canonical morphism $\sigma \colon \mr{Q}^n_{C}(\mcF) \to \Sym^n(C)$ that maps any $k$-rational point $[q: \mcF \to \mcH]$ to the effective divisor
\[
\sdiv (\mcH) := \sum_{p\in C}  \dim_k (\mcH_p) \, p \, ,
\]
where $\mcH_p$ is the stalk of $\mcH$ at $p$. 
\end{prop}

\begin{rmk}\label{ressigma}
If $U\subseteq C$ is open, then $\sigma^{-1}(\Sym^n(U))$ can be naturally identified with $\mr{Q}^n_U(\mcE)$. 
Hereafter, the morphism induced by $\sigma$ will be denoted by $\sigma_U \colon \mr{Q}^n_U(\mcE) \to \Sym^n(U)$. 
\end{rmk}

When $\mcE = \mcO_C$, $\mr{Q}^n_U(\mcO_C) = \Hilb_{U/k}^n$ and this morphism is the Hilbert-Chow morphism $\rho_U: \Hilb_{U/k}^n \to \Sym^n(U)$ of \cite{FG05}. Notice that $\rho_U$ is an isomorphism. 
Therefore, for any $\mcE$, $\sigma_U$ factors through $\rho_U$ via the morphism $\tau_U := \rho_U^{-1} \circ \sigma_U: \mr{Q}^n_U(\mcE) \to \Hilb_{U/k}^n$. 


\subsection{Explicit construction of $\tau_U$} \label{explsigma}
For later use, we provide here an explicit construction of $\tau_U$.

Let us consider the universal quotient $u_0: \mcE_{\mr{Q}^n_U(\mcE)} \to \mcQ_0$ associated to $\mr{Q}^n_U(\mcE)$. If $\mcK = \ker(u_0)$ and $\iota: \mcK \to \mcE_{\mr{Q}^n_U(\mcE)}$ is the inclusion, then we have the short exact sequence
\[
0 \to \mcK \xto{\iota} \mcE_{\mr{Q}^n_U(\mcE)} \xto{u_0} \mcQ_0 \to 0 \, .
\]
Since both $\mcE_{\mr{Q}^n_U(\mcE)}$ and $\mcQ_0$ are flat over $\mr{Q}^n_U(\mcE)$, $\mcK$ is flat over $\mr{Q}^n_U(\mcE)$ too.
Moreover, the restriction of $\mcK$ to $U_q$ is locally free, for all $q \in \mr{Q}^n_U(\mcE)$. 
It follows that $\mcK$ is a locally free sheaf of rank $r=\rk(\mcE)$. 

Let $\wedge^r (\iota): \wedge^r (\mcK) \to \wedge^r(\mcE_{\mr{Q}^n_U(\mcE)})$ be the $r$-th exterior power of $\iota$.
By tensoring it with $\wedge^r(\mcE_{\mr{Q}^n_U(\mcE)})^\vee$,  we get a short exact sequence
\[
0 \to \wedge^r (\mcK) \otimes \wedge^r(\mcE_{\mr{Q}^n_U(\mcE)})^\vee \to \mcO_{\mr{Q}^n_U(\mcE) \times U} \to \mcG \to 0 \, .
\]
Notice that $\mcG$ is flat over $\mr{Q}^n_U(\mcE)$, since $\mcO_{\mr{Q}^n_U(\mcE) \times U}$ is $\mr{Q}^n_U(\mcE)$-flat and
$\wedge^r(\iota)$ remains injective when restricted to every fiber (see \cite{Mat80}, Thm.\ 49 and its corollaries).
Moreover, the Hilbert polynomial of the restriction of $\mcG$ to every fiber is equal to $n$; indeed, the elementary divisor theorem for PIDs implies that the restriction of $\mcG$ to every fiber is isomorphic to the restriction of $\mcQ_0$ to the same fiber. 
Therefore, the quotient $\mcO_{\mr{Q}^n_U(\mcE) \times U} \to \mcG$ corresponds to a morphism $\mr{Q}^n_U(\mcE) \to \Hilb_{U/k}^n$, which is exactly the morphism $\tau_U$ defined in Section \ref{sigma}.



\subsection{The case $n=1$}\label{n=1}
The following result should be well known, but we include it here for lack of a suitable reference.

\begin{prop} \label{pbund}
The Quot scheme $\mr{Q}^1_C(\mcE)$ is isomorphic to the projective space bundle $\pr(\mcE)$.
\end{prop}

\proof
In order to simplify the notation, let us denote $\mr{Q}^1_C(\mcE)$ by $Q$, and the universal quotient over $Q \times C$ by $u: \mcE_Q \to \mcQ$.

Let $\sigma: Q \to C$ be the morphism introduced in Section \ref{sigma} and let
$\phi = (\id_Q, \sigma): Q \to Q \times C$ be the morphism with components the identity of $Q$ and $\sigma$, respectively. Then the pullback of $u$ via $\phi$ gives a quotient 
$\phi^\ast u: \phi^\ast\mcE_Q \to \phi^\ast\mcQ$. Notice that $\phi^\ast\mcQ \otimes {k}(q) \cong k(q)$ for any $q \in Q$ (where $k(q)$ is the residue field of $q$), therefore $\phi^\ast\mcQ$ is locally free of rank $1$.
Since $\phi^\ast\mcE_Q = \sigma^\ast\mcE$, by \cite[Prop.\ 7.12]{Har77} we obtain a morphism $\psi: Q \to \pr(\mcE)$.

On the other hand, let $\pi \colon \pr (\mcE) \to C$ be the  projection associated to $\mcE$, and let $\pi^*\mcE \to \mcL$ be the universal quotient over $\pr(\mcE)$. Let us consider the subscheme $\tilde{\Delta} \subset \pr (\mcE) \times C$ whose structure sheaf is $(\pi \times \id_C)^\ast \mcO_\Delta$, where $\pi \times \id_C: \pr(\mcE) \times C \to C \times C$ is the morphism with components $\pi$ and the identity of $C$ respectively, and $\Delta \subset C \times C$ is the diagonal.
The tensor product of the quotient $\mcO_{\pr(\mcE)\times C} \to \mcO_{\tilde{\Delta}}$ with $\mcE_{\pr(\mcE)}$ gives a surjection $q: \mcE_{\pr(\mcE)} \to \mcE_{\pr(\mcE)} \otimes \mcO_{\tilde{\Delta}}$.

If $\mr{pr}_1: \pr(\mcE) \times C \to \pr(\mcE)$ is the projection, then there is an isomorphism $\mcE_{\pr(\mcE)} \otimes \mcO_{\tilde{\Delta}} \cong (\mr{pr}_1)^\ast \pi^\ast \mcE \otimes \mcO_{\tilde{\Delta}}$.
Therefore, we can compose $q$ with the morphism $(\mr{pr}_1)^\ast \pi^\ast \mcE \otimes \mcO_{\tilde{\Delta}} \to (\mr{pr}_1)^\ast\mcL 
\otimes \mcO_{\tilde{\Delta}}$, and we obtain a surjective morphism $\mcE_{\pr(\mcE)} \to (\mr{pr}_1)^\ast\mcL \otimes \mcO_{\tilde{\Delta}}$.
Since $(\mr{pr}_1)^\ast\mcL \otimes \mcO_{\tilde{\Delta}}$ is flat over $\pr(\mcE)$ and has constant Hilbert polynomial $1$, this surjection corresponds to a morphism $\pr(\mcE) \to Q$, which is the inverse of $\psi$ by construction. 
\endproof

We conclude this section with the following result, from which we deduce that in general $\mr{Q}^n_C(\mcE)$
depends on $\mcE$, if $\rk(\mcE)\geq 2$.

\begin{prop}
Let $\mcE$ and $\mcE'$ be two locally free coherent $\mcO_C$-modules
of the same rank $r\geq 2$. Assume that one of the following conditions holds true:
\begin{enumerate}[label=(\roman*)]
\item
the genus of $C$ is greater than or equal to $1$;
\item
$r>2$.
\end{enumerate}
Then $\mr{Q}^1_C(\mcE) \cong \mr{Q}^1_C(\mcE')$ if and only if there exists an automorphism $\psi$ of $C$ and an invertible sheaf $\mcL$ on $C$, such that $\mcE' \cong \psi^\ast\mcE \otimes \mcL$. 
\end{prop}

\proof
First, assume that $\mr{Q}^1_C(\mcE) \cong \mr{Q}^1_C(\mcE')$. By Proposition \ref{pbund}, we thus have an isomorphism $\phi: \pr(\mcE') \iso \pr(\mcE)$. Let us denote the projections of these bundles by $\pi: \pr(\mcE) \to C$ and $\pi': \pr(\mcE') \to C$. Under the hypothesis (i) or (ii), the morphism $\pi \circ \varphi: \pr(\mcE') \to C$ is constant on the fibers of $\pi'$. Therefore, there exists an automorphism $\psi$ of $C$ such that the diagram
\[
\begin{tikzcd}
	\pr(\mcE') \arrow{r}{\phi} \arrow{d}{\pi'} & \pr(\mcE) \arrow{d}{\pi} \\
	C \arrow{r}{\psi} & C
\end{tikzcd}
\]
commutes. Since there is also a cartesian diagram
\begin{equation} \label{pbuncart}
\begin{tikzcd}
	\pr(g^\ast\mcE) \arrow{r}{\cong} \arrow{d}{p} & \pr(\mcE) \arrow{d}{\pi} \\
	C \arrow{r}{\psi} & C
\end{tikzcd}
\end{equation}
where $p$ is the canonical projection, it follows that $\pi': \pr(\mcE') \to C$ and $p: \pr(g^\ast\mcE) \to C$ are isomorphic over $C$. By \cite[\S 2.7, Exer.\ 7.9]{Har77}, we deduce that $\mcE' \cong \psi^\ast\mcE \otimes \mcL$ for some invertible sheaf $\mcL$ on $C$.

The inverse implication directly follows from \cite[\S 2.7, Exer.\ 7.9]{Har77} and the diagram \eqref{pbuncart}.
\endproof


\section{The fibers of $\sigma$} \label{sigmafiber}
In this section we describe the fibers of the morphism $\sigma$ introduced in Section \ref{sigma}.
Throughout the section, $C$ denotes a smooth projective curve over $k$ and $\mcE$ is a coherent locally free $\mcO_C$-module of rank $r$.

\begin{prop} \label{samefiber}
The fiber of $\sigma \colon \mr{Q}^n_C(\mcE) \to \Sym^n ( C )$ over a point $D\in \Sym^n ( C )$ is isomorphic 
to the fiber of  the analogous morphism $\mr{Q}^n_C(\mcO_C^{\oplus r}) \to \Sym^n ( C )$ over the same point.
\end{prop}

\proof
From Remark \ref{ressigma} we have that $\sigma^{-1}(\Sym^n(U))$ depends only on $\mcE|_U$, for any $U\subseteq C$ open. Then the proposition follows from the fact that for any $D\in \Sym^n(C)$, there exists an open subset $U\subseteq C$ such that $D \in \Sym^n(U)$ and $\mcE|_U$ is trivial. 

In order to see this, let $V$ be an open affine subset of $C$, such that
$D \in \Sym^n(V)$. Then, by \cite[Thm.\ 1]{Ser58}, $\mcE|_V \cong \mcO_V^{\oplus (r-1)} \oplus \mcL$, where $\mcL$ is an invertible $\mcO_V$-module. 
Let us consider the short exact sequence $0 \to \mcL(-D) \to \mcL \to \mcL \otimes \mcO_D \to 0$, and the associated exact cohomology sequence, $0 \to \h^0(V, \mcL(-D)) \to \h^0(V, \mcL) \to \h^0(V, \mcL \otimes \mcO_D) \to 0$.
We deduce that there exists $s \in \h^0(V,\mcL)$ such that $s(x) \neq 0$, for all $x\in \Supp(D)$. Hence $\mcL$ (and consequently $\mcE$) is trivial on the open set $U = V \setminus \Supp(s)$.
\endproof

\begin{dfn} \label{fibsigma}
For any $p\in C$, let us define $F_{n,r}(p) := \sigma^{-1}(np)$.
More generally, for any $D\in \Sym^n(C)$, we define $F_r(D):= \sigma^{-1}(D)$.
\end{dfn}

\begin{prop}
Let $D= a_1p_1 + \ldots + a_mp_m \in \Sym^n(C)$, with $p_1, \dots, p_m$ pairwise distinct. Then 
\[
F_r(D) \cong F_{a_1,r}(p_1) \times \dots \times F_{a_m,r}(p_m) \, .
\]
\end{prop}

\proof
There is a natural morphism $F_r(D) \to F_{a_1,r}(p_1) \times \dots \times F_{a_m,r}(p_m)$, which is defined in the following way on $k$-rational points.
For any quotient $[q: \mcE \to \mcH]$ in $F_r(D)$, we have a splitting 
$\mcH = \oplus_{i=1}^m \mcH_{p_i}$, where $\mcH_{p_i}$ is a skyscraper sheaf on $C$, which is supported in $\{p_i\}$. Therefore $q=\oplus_{i=1}^m q_i$, with $q_i: \mcE \to \mcH_{p_i}$. Then $[q: \mcE \to \mcH]$ is mapped to $[q_1: \mcE \to \mcH_{p_1}] \times \dots \times [q_m: \mcE \to \mcH_{p_m}]$. Clearly this is an isomorphism. 
\endproof


\section{The class of $\mr{Q}^n_C(\mcE)$ in $\Kvar$} \label{coresect} 

In this section we prove our main theorem. 

\begin{thm}\label{mainthm}
Let $C$ be a smooth projective curve over $k$.
Let $\mcE$ be a coherent locally free $\mcO_C$-module of rank $r$. Then, for any non-negative integer $n$, the equality 
\[
[\mr{Q}^n_C(\mcE)] = [\mr{Q}^n_C(\mcO_C^{\oplus r})] 
\] 
holds true in the Grothendieck group $\Kvar$ of $k$-varieties. 
\end{thm}

The proof will be divided into several steps.

\subsubsection*{Step 1} 
In order to make the proof clearer, we first fix our notation (see also Section \ref{notations}). Let $U\subseteq C$ be a fixed open subset such that $\mcE|_{U} \cong \mcO_U^{\oplus r}$, and let 
$C\setminus U = \{ p_1, \ldots , p_N \}$. Then $\Sym^n ( C )$ is the set-theoretic disjoint union of the locally closed subsets
\[
Z_\mb{a}:= \{ E \in \Sym^n ( C ) \mid \Supp ( E - a_1p_1 -\ldots  - a_Np_N ) \subset U \} \,,
\]
for $\mb{a} \in A := \{(a_1, \ldots , a_N)\in \N^N \mid a_1 + \ldots + a_N \leq n \}$. 
Notice that 
\[
Z_\mb{a} \cong \Sym^{n-|\mb{a}|} ( U ) \,,
\]
where $|\mb{a}| := a_1 + \dots + a_N$.
 
For any $\mb{a} \in A$, we denote by $\mr{Q}_\mb{a}(\mcE)$ the preimage of $Z_\mb{a}$ under the morphism $\sigma: \mr{Q}^n_C(\mcE) \to \Sym^n(C)$ of Section \ref{sigma}, with the reduced 
subscheme structure. 

\begin{rmk} \label{strat}
Using the relations in the Grothendieck group of varieties, the decomposition of $\mr{Q}^n_C(\mcE)$ into its locally closed subsets $\mr{Q}_\mb{a}(\mcE)$ yields the equality
\[
[\mr{Q}^n_C(\mcE)] = \sum_{\mb{a} \in A} [\mr{Q}_\mb{a}(\mcE)] 
\]
in $\Kvar$.
\end{rmk}

Finally, we denote by $D$ the divisor $a_1p_1 + \dots + a_Np_N \in \Sym^{|\mb{a}|}(C)$ corresponding to $\mb{a} \in A$. Associated to this effective divisor we have the fiber $F_r(D) \subset \mr{Q}^{|\mb{a}|}_C(\mcE)$, as in Definition \ref{fibsigma}.

\subsubsection*{Step 2}
The core of our proof is the following proposition.

\begin{prop} \label{stratiso}
For any $\mb{a} \in A$, there is a natural isomorphism 
\[
\mr{Q}_\mb{a}(\mcE) \iso \mr{Q}^{n-|\mb{a}|}_U(\mcE) \times F_{r}(D) \, .
\]
\end{prop}

The idea behind this proposition is that any quotient of $\mcE$ in $\mr{Q}_\mb{a}(\mcE)$
can be obtained by glueing a quotient supported in $U$ and a quotient supported on $\{ p_1, \ldots , p_N\}$.

Before proving Proposition \eqref{stratiso} in Step 3, we need the following result in order to define the morphism $\mr{Q}_\mb{a}(\mcE) \to \mr{Q}^{n-|\mb{a}|}_U(\mcE)$.

\begin{lmm}\label{suppQ|proper}
Let $\mcQ$ be the universal quotient associated to $\mr{Q}^n_C(\mcE)$, and let $i: U\times \mr{Q}_\mb{a}(\mcE) \hookrightarrow C \times \mr{Q}^n_C(\mcE)$ be the inclusion.
For any $\mb{a} \in A$, the support $\Supp(i^\ast \mcQ)$ is proper over $\mr{Q}_\mb{a}(\mcE)$.
\end{lmm}

\proof 
We apply the valuative criterion of properness to the restriction $\phi: \Supp(i^\ast \mcQ) \to \mr{Q}_\mb{a}(\mcE)$ of the projection $U \times \mr{Q}_\mb{a}(\mcE) \to \mr{Q}_\mb{a}(\mcE)$.
So, let $R$ be a valuation ring, and let $K$ be its  quotient field. Assume we are given a commutative diagram
\[
\begin{tikzcd}
\Spec(K) \arrow{d} \arrow{r} & \Supp(i^\ast\mcQ) \arrow{d}{\phi} \\
\Spec(R) \arrow{r} & \mr{Q}_\mb{a}(\mcE)
\end{tikzcd}
\] 
where $\Spec(K) \to \Spec(R)$ is the morphism induced by the inclusion $R \into K$.

Since $\Supp (i^\ast \mcQ ) \subset C \times \mr{Q}_\mb{a}(\mcE)$ and $C \times \mr{Q}_\mb{a}(\mcE)$ is proper over $\mr{Q}_\mb{a}(\mcE)$, there is a unique morphism $\psi: \Spec(R) \to C \times \mr{Q}_\mb{a}(\mcE)$ such that the diagram 
\[
\begin{tikzcd}
\Spec(K) \arrow{d} \arrow{r} & \Supp (i^\ast \mcQ) \arrow{d}[near start]{\phi} \arrow[hookrightarrow]{r} & C \times \mr{Q}_\mb{a}(\mcE) \\
\Spec(R) \arrow{r} \arrow[dotted]{rru}[near start]{\psi} &  \mr{Q}_\mb{a}(\mcE) \arrow[leftarrow]{ru} &
\end{tikzcd}
\]
commutes. The claim follows, if we prove that the image of $\psi$ is contained in $U \times \mr{Q}_\mb{a}(\mcE)$.

Let us consider the composition
\[
\psi': \Spec(R) \xrightarrow{\psi} C \times \mr{Q}_\mb{a}(\mcE) \xrightarrow{\id_C \times (\tau|_{\mr{Q}_\mb{a}(\mcE)})} C \times \rho^{-1}(Z_\mb{a}) \,,
\]
where $\tau: \mr{Q}^n_C(\mcE) \to \Hilb^n_{C/k}$ is the morphism defined in Section \ref{sigma}.
The universal ideal sheaf of $\Hilb_{C/k}^n$ restricted to  $C \times \rho^{-1}(Z_\mb{a})$ is of the form $\mcO_{C\times \rho^{-1}(Z_\mb{a})}(-\mcD - \mcD')$, for $\mcD=a_1\{p_1\}\times \rho^{-1}(Z_\mb{a}) +\ldots + a_N\{p_N\}\times \rho^{-1}(Z_\mb{a})$ and $\mcD'$ an effective Weil divisor (notice that $C\times \rho^{-1}(Z_\mb{a})$
is smooth), such that $\mcD \cap \mcD' = \emptyset$. 
Since the image of the composition
\[
\Spec(K) \to \Supp(i^\ast \mcQ) \into C \times \mr{Q}_\mb{a}(\mcE) \xrightarrow{\id_C \times (\tau|_{\mr{Q}_\mb{a}(\mcE)})} C \times \rho^{-1}(Z_\mb{a})
\]
is contained in $\Supp(\mcD') \subset U \times \rho^{-1}(Z_\mb{a})$, the same holds for the image of $\psi'$. Therefore, the image of $\psi$ lies in $U \times \mr{Q}_\mb{a}(\mcE)$, as claimed.
\endproof

\subsubsection*{Step 3}

Using Lemma \ref{suppQ|proper}, we can now prove Proposition \ref{stratiso}.

\proof[Proof of Proposition \ref{stratiso}]
As above, let $u: \mcE_{\mr{Q}^n_C(\mcE)} \to \mcQ$ be the universal quotient associated to $\mr{Q}^n_C(\mcE)$, and let $i: U\times \mr{Q}_\mb{a}(\mcE) \into C \times \mr{Q}^n_C(\mcE)$ be the inclusion.
By Lemma \ref{suppQ|proper}, the quotient $i^\ast u: i^\ast(\mcE_{\mr{Q}^n_C(\mcE)}) \to i^\ast\mcQ$ yields a natural morphism $f_0: \mr{Q}_\mb{a}(\mcE) \to \mr{Q}^{n-|\mb{a}|}_U(\mcE)$. 

In order to define a morphism $f_\infty: \mr{Q}_\mb{a}(\mcE) \to F_{r}(D)$, let us consider the open neighbourhood $(C \times \mr{Q}_\mb{a}(\mcE)) \setminus \Supp(i^\ast\mcQ)$ of $\cup_{\lambda=1}^N \{ p_\lambda \} \times \mr{Q}_\mb{a}(\mcE)$, together with its inclusion $j$ into $C \times \mr{Q}_\mb{a}(\mcE)$.
By composing the pullback of $u$ to $\mr{Q}_\mb{a}(\mcE)$ with the unit of the adjunction $j^\ast \dashv j_\ast$, we get a surjective morphism
\[
q: \mcE_{\mr{Q}_\mb{a}(\mcE)} \to \mcQ_{\mr{Q}_\mb{a}(\mcE)} \to j_\ast j^\ast(\mcQ_{\mr{Q}_\mb{a}(\mcE)})
\]
of coherent $\mcO_{C \times \mr{Q}_\mb{a}(\mcE)}$-modules.
Notice that $j_\ast j^\ast(\mcQ_{\mr{Q}_\mb{a}(\mcE)})$ has constant Hilbert polynomial $|\mb{a}|$, hence it is flat over $\mr{Q}_\mb{a}(\mcE)$. 
Therefore, $q$ is associated to a natural morphism $\mr{Q}_\mb{a}(\mcE) \to \mr{Q}^{|\mb{a}|}_C(\mcE)$, whose image is contained in $F_r(D)$. Thus we obtain a morphism $f_{\infty}: \mr{Q}_\mb{a}(\mcE) \to F_r(D)$, and the morphism in \eqref{stratiso} is $(f_0, f_\infty)$.

To prove that $(f_0, f_\infty)$ is an isomorphism, we exhibit its inverse, as follows.
Let $u_0: \mcE_{\mr{Q}^{n-|\mb{a}|}_U(\mcE)} \to \mcQ_0$ be the universal quotient of $\mr{Q}^{n-|\mb{a}|}_U(\mcE)$, and let $u_\infty: \mcE_{F_r(D)} \to \mcQ_\infty$ be the pullback of the universal quotient of $\mr{Q}^{|\mb{a}|}_C(\mcE)$ to $C \times F_r(D) \subset C \times \mr{Q}^{|\mb{a}|}_C(\mcE)$. 
In the following, we view $u_0$ as a family of quotients of $\mcE$ supported in $U$ (see Lemma \ref{GN}).
Let us denote the projections by $\mr{pr}_{12}: C \times \mr{Q}^{n-|\mb{a}|}_U(\mcE) \times F_r(D) \to 
C \times \mr{Q}^{n-|\mb{a}|}_U(\mcE)$ and $\mr{pr}_{13}: C \times \mr{Q}^{n-|\mb{a}|}_U(\mcE) \times F_r(D) \to C  \times F_r(D)$. 
Then 
\[
(\mr{pr}_{12})^\ast u_0 \oplus (\mr{pr}_{13})^\ast u_\infty: \mcE_{\mr{Q}^{n-|\mb{a}|}_U(\mcE) \times F_r(D)} \to (\mr{pr}_{12})^\ast \mcQ_0 \oplus (\mr{pr}_{13})^\ast \mcQ_\infty
\]
is a family of quotients of $\mcE$, parametrized by $\mr{Q}^{n-|\mb{a}|}_U(\mcE) \times F_r(D)$, with constant Hilbert polynomial equal to $n$. The associated morphism 
$\mr{Q}^{n-|\mb{a}|}_U(\mcE) \times F_r(D) \to \mr{Q}_C^n(\mcE)$ is the inverse morphism of $(f_0, f_\infty)$.
\endproof

\subsubsection*{Step 4}
We can finally conclude the proof of Theorem \ref{mainthm}.

By Remark \ref{strat} and Proposition \ref{stratiso}, we have the following equalities in $\Kvar$:
\[
[\mr{Q}^n_C(\mcE)] = \sum_{\mb{a} \in A} [\mr{Q}_\mb{a}(\mcE)] = \sum_{\mb{a} \in A} [\mr{Q}^{n-\mb{a}}_U(\mcE)] [F_r(D)] \,.
\]
In particular, this is true also for $\mcE = \mcO_C^{\oplus r}$.

Now, $\mcE|_U$ is trivial, therefore $[\mr{Q}^{n-\mb{a}}_U(\mcE)] = [\mr{Q}^{n-\mb{a}}_U(\mcO_C^{\oplus r})]$. Theorem \ref{mainthm} thus follows from Proposition \ref{samefiber}.


\subsection{Explicit computation}
We provide an explicit formula for the class $[\mr{Q}^n_C(\mcE)] \in \Kvar$
in terms of the classes $[\Sym^m(C)]$.

\begin{prop}\label{class}
For any non-negative integer $n$, the equality 
\[
[\mr{Q}^n_C(\mcE)] = \sum_{\mb{n} \in \N^r , \, |\mb{n}|=n} [\Sym^{n_1}(C)]\cdot \ldots \cdot [\Sym^{n_r}(C)] \cdot [\mathbb{A}^1_k]^{d_\mb{n}}
\] 
holds true in $\Kvar$, where $d_\mb{n} := \sum_{i=1}^r (i-1)n_i$.
\end{prop}

\proof
From Theorem \ref{mainthm} we have that $[\mr{Q}^n_C(\mcE)] = [\mr{Q}^n_C(\mcO_C^{\oplus r})]$. 
Then the result follows directly from \cite{Bif89}.
\endproof

\begin{rmk}
From the previous formula we can determine the Poincar\'e polynomial 
of $\mr{Q}^n_C(\mcE)$ (for $\ell$-adic cohomology, where $\ell \neq \car(k)$ is a prime) as follows (cf.\ also \cite{BGL94}). By \cite{Mac62}, the Poincar\'e polynomial $P(\Sym^m(C); t)$ of $\Sym^m(C)$ is the coefficient of $u^m$ in the expansion of 
\[
\frac{(1+tu)^{2g}}{(1-u)(1-t^2u)} \, ,
\]
where $g$ is the genus of $C$.
Then, for $E(t,u):= \sum_{n=0}^\infty P(\mr{Q}^n_C(\mcE);t)\, u^n$, we have:
\begin{align*}
E(t,u) &= \sum_{n=0}^\infty \sum_{\mb{n} \in \N^r , \, |\mb{n}|=n} P(\Sym^{n_1}(C);t) \cdots P(\Sym^{n_r}(C);t)\, t^{2d_\mb{n}}u^n \\
&= \sum_{\mb{n} \in \N^r} \prod_{i=1}^{r} P(\Sym^{n_i}(C);t)\, u^{n_i} t^{2(i-1)n_i} \\
&= \prod_{i=0}^{r-1} \frac{(1+t^{2i+1}u)^{2g}}{(1-t^{2i}u)(1-t^{2i+2}u)} \,.
\end{align*} 
\end{rmk}

\begin{cor}\label{irreducibility}
The Quot scheme $\mr{Q}^n_C(\mcE)$ is irreducible.
\end{cor}
\proof
Since $\mr{Q}^n_C(\mcE)$ is smooth, it suffices to show that the coefficient of $t^{0}$
in the Poincar\'e polynomial of $\mr{Q}^n_C(\mcE)$ is $1$. To this aim, notice that, for any $\mb{n} \in \N^r$ with $|\mb{n}|=n$, we have that $d_\mb{n} =\sum_{i=1}^r (i-1)n_i =0$ only for $\mb{n}=(n,0,\ldots, 0)$. 
The claim now follows from Proposition \ref{class}. 
\endproof

\bibliographystyle{amsalpha}
\bibliography{bib}

\end{document}